\documentclass[10pt]{amsart}
\usepackage[dvips,ps,color,all]{xy}
\usepackage{graphicx,amssymb,latexsym,amsfonts,verbatim,amscd,amsmath,amsthm}
\usepackage[colorlinks=true,linkcolor=black]{hyperref}
\providecommand{\dis}{\displaystyle}

\theoremstyle{plain}
\newtheorem{thm}{Theorem}

\newtheorem{cor}[thm]{Corollary}
\newtheorem{lemma}[thm]{Lemma}
\newtheorem{prop}[thm]{Proposition}
\newtheorem{defi}[thm]{Definition}
\theoremstyle{remark}
\newtheorem{rem}{Remark}
\newtheorem{rems}{Remarks}
\newtheorem{ex}{Example}

\begin{document}
\title{On special $p$-Borel fixed ideals}
\author{Achilleas Sinefakopoulos}
\date{February 20, 2007}

\begin{abstract}
 We define the reduced horseshoe resolution and the notion of
conjoined pairs of ideals in order to study the minimal graded
free resolution of a class of $p$-Borel ideals and recover
Pardue's regularity formula for them. It will follow from our
technique that the graded betti numbers of these ideals do not
depend on the characteristic of the base field $\Bbbk$.
\end{abstract}
\maketitle
\section{Introduction} \label{S:I}

The study of $p$-Borel fixed ideals is a very interesting and
 fascinating problem. One could safely argue that in characteristic
 $p$ very few results are known, in contrast to the case of
 characteristic zero, where we can describe a minimal
 graded free resolution of any Borel fixed ideal, determine its
 regularity, find its graded Betti numbers and more.

It was conjectured in \cite{Par} that the regularity $reg(I)$ of a
principal $p$-Borel ideal $I$ is equal to the maximum of some
numbers given by a rather complicated formula. In \cite{Ara} it
was proved that $reg(I)$ is larger than or equal to this maximum,
while in \cite{HerPop} the authors prove the opposite inequality
(see also \cite{HerPopV}).

Another known result is the computation in \cite{Ene} of the
Koszul homology of some \emph{special} $p$-Borel ideals, which we
shall define below, while a more recent result is their
CW-resolution given in \cite{Wel} using algebraic discrete Morse
theory. In both papers there are proofs of a formula that gives
their regularity, which agrees with Pardue's formula for principal
$p$-Borel ideals.

Here we show how one can use the horseshoe lemma to get the form
of the minimal graded free resolution of these \emph{special}
$p$-Borel ideals in an elementary way. Furthermore, we verify
Pardue's regularity formula at the same time. Our idea was born
from the observation of the Betti diagrams of several such ideals
in MACAULAY 2 \cite{Mike}.

This paper is organized as follows:

 In  Section 2, we introduce the \emph{reduced
horseshoe resolution}, which will help us deduce the minimal
graded free resolution of $R/IJ$ from the minimal resolutions of
$R/I$ and $R/J$, when the pair $(I,J)$ of ideals in $R$ satisfies
certain properties.

In Section 3, we study the minimal graded free resolution over $R$
of a class of ideals, which we call \emph{special}.
\[
    I=\prod_{j=1}^{s} I_{j}^{[p_j]},
\]
where
\[
I_j=(x_1,x_2,...,x_{\ell_j})^{a_j} \quad \text{and} \quad
I_j^{[p_j]}=(x_1^{p_j},x_2^{p_j},...,x_{\ell_j}^{p_j})^{a_j}
 \]
 with
 \[
   n=\ell_1 \geq \ell_2\geq ...\geq \ell_s \geq 1 \quad \text{and} \quad 0\leq a_j
   <\dfrac{p_{j+1}}{p_j},
    \]
where the numbers $\dfrac{p_{j+1}}{p_j}$ are integers $>1$ for
$j=1,...,s$.

We call such ideals \emph{special}. In particular, if
$p_j=p^{r_j}$ for $j=1,2,...,s$, where $r_s>...>r_2>r_1\geq 0$ and
$p$ is prime, we call them \emph{special} $p$-Borel ideals.

 In Section 4, we construct a polyhedral cell
    complex that supports a minimal free resolution of some special
    ($p$-Borel) ideals.

 Finally, in section 5, we examine the iterated mapping cone
    construction.

\section{Reduced Horseshoe Resolution and Conjoined Pairs of Ideals}\label{RHR:HC}

    All ideals in this paper are considered to be \emph{monomial}
ideals. We work over the polynomial ring
$R=\Bbbk[x_1,x_2,...,x_n]$. For small $n$ we may use the letters
$a,b,c,d,...$ instead of $x_1,x_2,x_3,x_4,...$.

Let $A \subset B$ be two ideals in $R$ and assume that the minimal
graded free resolutions of $A/B$ and $A$ are of the form
\begin{small}
\[
 \xymatrix{
0\ar[r]&G_m \ar[r]^{d_m^{''}}&\dots\ar[r]& G_2\ar[r]^{d_2^{'}} &
G_1\ar[r]^{d_1^{'}}&F_1\ar[r]^{\epsilon_{0}^{'}}& A/B \ar[r]&0
                      }
                      \]
\end{small}
and
\begin{small}
\[
 \xymatrix{
0\ar[r]&F_n \ar[r]^{d_n^{''}}&\dots\ar[r]& F_2\ar[r]^{d_2^{''}} &
F_1\ar[r]^{d_1^{''}} & R\ar[r]^{\epsilon_{0}^{''}}& R/A \ar[r]&0
                      }
                      \]
\end{small}

Let $\dis F_i=\oplus_{i=1}^{b_k} R(-\alpha_{i,k})$ for
$i=1,2,...,n$ and $k=1,2,...,b_i$. and set
\[
\begin{aligned}\notag
   m_i(A)&= \text{min}\{\alpha_{i,k}| k=1,2,..,b_i \}\\\notag
   M_i(A)&= \text{max}\{\alpha_{i,k}| k=1,2,..,b_i \}
\end{aligned}
 \]
for $i=1,2,...,n$. Then the horseshoe lemma associated with the
following short exact sequence
\begin{small}
\[
 \xymatrix{
           0\ar[r]& A/B \ar[r]^{\psi} & R/B \ar[r]^{\phi} &R/A\ar[r]&0
                      }
                      \]
                      \end{small}
gives us a free resolution of $R/B$,
\begin{small}
\[
 \xymatrix{
                               &0\ar[d]                                      &0\ar[d]                          &0\ar[d]                      &0\ar[d]                          &0\ar[d]         &\\
                   \dots\ar[r] & G_3 \ar[d]^{\psi_3}\ar[r]^{d_3^{'}}                  & G_2\ar[d]^{\psi_2}\ar[r]^{d_2^{'}}       & G_1 \ar[d]^{\psi_1}\ar[r]^{d_1^{'}}          & F_1\ar[d]^{\psi_0}\ar[r]^{\epsilon_{0}^{'}}      & A/B \ar[d]^{\psi}\ar[r]&0\\
                   \dots\ar[r] & G_3 \oplus F_3 \ar[d]^{\phi_3} \ar[r]^{d_3}                &G_2 \oplus F_2 \ar[d]^{\phi_2}\ar[r]^{d_2}      & G_1 \oplus F_1 \ar[d]^{\phi_1}\ar[r]^{d_1}  &  F_1 \oplus R \ar[d]^{\phi_0}\ar[r]^{\epsilon_0}       & R/B \ar[d]^{\phi}\ar[r]&0\\
                 \dots\ar[r] & F_3\ar@{-->}[uur] \ar[d]\ar[r]^{d_3^{''}}                  & F_2\ar@{-->}[uur]\ar[d]\ar[r]^{d_2^{''}}       & F_1\ar@{-->}[uur]\ar[d]\ar[r]^{d_1^{''}}   & R\ar[ur]^{\pi}\ar[d]\ar[r]^{\epsilon_{0}^{''}}  & R/A \ar[d]\ar[r]&0\\
                             & 0                                             &0                                 &0                              &0                                  &0             & }
                      \]
                      \end{small}
\newline
The differential map  $\epsilon_0 : F_1 \oplus R \to R/B$ is
defined by $\epsilon_0 (x,y)=\psi \epsilon_0^{'} (x)+\pi(y)$ for
$x \in F_1$ and $y \in R$, while the maps
 $d_k: G_k \oplus F_k \to G_{k-1} \oplus F_{k-1}$ for $k>1$ are given by the following matrix
                      \[
                        d_k=\begin{bmatrix}
                        d_k^{'} &\lambda_k  \\
                        0 & d_k^{''} \\
                        \end{bmatrix},\]
where the maps $\lambda_k$ are the ones denoted by the dashed
arrows in the above commutative diagram. Moreover, the above maps
must satisfy the following conditions (see, e.g. \cite{Car},
p.79-80)
\begin{align}\notag
 \epsilon_0^{''}&=\phi \pi, \\\notag
 \psi \epsilon_{0}^{'} \lambda_1 +\pi d_{1}^{''}&=0 \qquad \text{and}\\\notag
d_{k-1}^{'} \lambda_{k} + \lambda_{k-1} d_{k}^{''}&=0 \qquad
\text{for} \qquad k>1.\notag
\end{align}
We may assume that $\epsilon_0^{'}:=\pi_{|A} d_1^{''}$, because
$Im(d_1^{''})=A$. In order to define $\lambda_1$, we choose a
basis element $e$ of $F_1$. Then $\pi
d_{1}^{''}(e)=d_{1}^{''}(e)+B$ and
\[
\psi \epsilon_{0}^{'} \lambda_1 (e)=\psi(\pi_{|A}
d_1^{''}(\lambda_1 (e)))=d_{1}^{''}(\lambda_1 (e))+B.
\]
 Thus, we need only make sure that $d_{1}^{''}(e)+d_{1}^{''}(\lambda_1 (e))$
is in $B$. Accordingly, we may define the map $\lambda_1 : F_1 \to
F_1$ such that $\lambda_1 (e)=-e$. Hence, the above horseshoe
resolution of $R/B$ is certainly not minimal. This leads us to
define the \emph{reduced horseshoe resolution} of $R/B$,
\begin{defi}
 Let $A$ and $B$ be two ideals in $R$ as above. Then the complex
\begin{small}
\[
 \xymatrix{
            \dots\ar[r]& G_3 \oplus F_3  \ar[rr]^{\begin{bmatrix}
                        d_3^{'} &\lambda_3  \\
                        0 & d_3^{''}
                        \end{bmatrix}}  &&G_2 \oplus F_2\ar[rr]^{\begin{bmatrix}
                        d_2^{'} &\lambda_2
                        \end{bmatrix}} &&G_1 \ar[rr]^{ d_1^{''}\lambda_1^{-1} d_1^{'}}  && R \ar[r]^{\pi} & R/B \ar[r]&0\\
 }
                      \]
                      \end{small}
is called the  \emph{reduced horseshoe resolution} of $R/B$ with
respect to $A$.
\end{defi}
It is easy to verify that this is a complex. Indeed, note that
$\pi d_1^{''} \lambda_1^{-1} d_1^{'}= -\psi
\epsilon_{0}^{'}d_1^{'} =0$, because $ \epsilon_{0}^{'}d_1^{'}=0$,
while
\[d_1^{''}\lambda_1^{-1} d_1^{'} \begin{bmatrix}
                        d_2^{'} &\lambda_2
                        \end{bmatrix} = \begin{bmatrix}
                       d_1^{''}\lambda_1^{-1}( d_1^{'} d_2^{'}) &d_1^{''} \lambda_1^{-1}(d_1^{'} \lambda_2)
                        \end{bmatrix}= \begin{bmatrix}
                       0 &-d_1^{''}d_2^{''}
                        \end{bmatrix}= \begin{bmatrix}
                       0 & 0
                        \end{bmatrix}.\]
The other relations follow immediately from the fact that $d_k
d_{k+1}=0$ for $k>1$.

\begin{rems}.
\begin{itemize}
    \item[(a)] Although the horseshoe resolution of $R/B$ is not
minimal, there is a chance that its reduced horseshoe resolution
is. For this to be true, the rank of the free $R$-module $G_1$
should be equal to the number of the minimal generators of $B$ and
the matrices that represent the maps $\lambda_k$ for $k>1$ should
not have any nonzero constant entry. We are already sure that the
matrices that represent the maps $d_k^{'}$ and $d_k^{''}$ do
satisfy the latter condition because we started with minimal
resolutions of $R/A$ and $A/B$.
    \item[(b)] If $A=I$ and $B=IJ$ for some ideal $I$ and $J$ in $R$,
we will talk about the \emph{reduced horseshoe resolution} of
$R/IJ$ with respect to the ordered pair $(I,J)$. More generally,
if $I_1$,....,$I_k$ $(k>1)$ are ideals in $R$, such that we know
the minimal graded free resolutions of $R/I_k$ and $I_k/ I_{k-1}
...I_{1}$, we will talk about the \emph{reduced horseshoe
resolution} of $R/(I_k ... I_2 I_1)$ with respect to the ordered
$k$-tuple $(I_k,...,I_2,I_1)$.
\end{itemize}
\end{rems}

Let $G(I)$ denote the unique minimal set of monomial generators of
a (monomial) ideal $I$.

\begin{defi} We call an ordered pair $(I,J)$ of monomial ideals
$I$ and $J$ in $R$ \emph{conjoined}, if the following conditions
are satisfied:
\begin{itemize}
   \item[(i)] $|G(IJ)|=|G(I)||G(J)|$.
\item[(ii)]There is a minimal presentation of $I$,
\begin{small}
 \[
\xymatrix{
   R^{s}\ar[r]^{\phi}& R^{t}\ar[r]^{\psi}&I\ar[r]&0
}
\]
\end{small}
such that all the entries of the matrix $\phi$ belong to $J$.
\end{itemize}
\end{defi}
\begin{ex}\label{EX:HC}
Let $I=(a^2 ,b^2, c^2)$ and $J=(a,b,c)$ in $R=\Bbbk[a,b,c]$. Then
$IJ=(a^3,a^2b,ab^2,b^3,a^2c,b^2c,ac^2,bc^2,c^3)$ and a minimal
presentation of $I$ is
\begin{small}
\[
\xymatrix{
   R\ar[rr]^{\begin{pmatrix}
     b^2 & -c^2 & 0 \\
     -a^2 & 0   & -c^2 \\
     0   & a^2 & b^2 \\
   \end{pmatrix}} & &R^2\ar[r]&(a^2,b^2, c^2)\ar[r]&0}.
\]
\end{small}
Hence the pair $(I,J)$ is \emph{conjoined}.
\end{ex}
The following lemma, which was inspired by lemma 2.2 of
\cite{Ene}, gives us a systematic way of constructing
\emph{conjoined} pairs of ideals.
\begin{lemma}\label{L1:HC}
Let $I$ and $J$ be two monomial ideals in
$R=\Bbbk[x_1,x_2,...,x_n]$ such that every element of $G(I)$ has
degree $d_1$ and every element of $G(J)$ has degree $d_2$. Suppose
that $(x_1^{k-1},....,x_m^{k-1}) \subset J$ for some integer $k>
1$, where $m=\text{max}\{i| x_i\quad \text{divides a minimal
generator of}\quad $I$\}$. Then the pair $(I^{[k]}, J)$ is
\emph{conjoined}, where $I^{[k]}$ is the ideal generated by $\{
u^k: u \in G(I)\}$.
\end{lemma}
\emph{Proof.} We first show that
\[
  G(I^{[k]}J)=\{u^{[k]}v:u \in G(I), v \in G(J)\}.
\]
Certainly, the set on the right hand side above is a generating
set of $I^{[k]}J$. We need to show that it is minimal. Assume that
there are $u \in G(I)$ and $ v \in G(J)$ such that $u^{[k]}v$  is
not a minimal generator. Then, we should have
\[
  u^{[k]}v={u'}^{[k]}v'w
\]
for some monomial $w$. This is impossible, because the above
relation together with $deg(u)=deg(u')$ and $deg(v)=deg(v')$
implies that $deg(w)=0$. Since $(x_1^{k-1},....,x_m^{k-1}) \subset
J$, the degree $d_2$ of every element in $G(J)$ is less than or
equal to $k-1$. Now if $\alpha_i$ and $\alpha_i'$ denote the
largest integers such that $x_i^{\alpha_i} $ divides $v$ and
$x_i^{\alpha_i^{'}}$ divides $v'$, then $0\leq |\alpha_i
-\alpha_i^{'}|\leq d_2 < k$ for each $i$. Since $k$ divides
$|\alpha_i -\alpha_i^{'}|$, we have $\alpha_i =\alpha_i^{'}$ for
each $i$ and so $v=v'$. Accordingly, we also have $u=u'$, and so
\[
  |G(I^{[k]}J)|=|G(I^{[k]})||G(J)|.
\]
Now consider a minimal presentation of $I$,
\begin{small}
 \[
\xymatrix{
   R^{s}\ar[r]^{(a_{i,j})}& R^{t}\ar[r]&I\ar[r]&0
}
\]
\end{small}
where all the entries $a_{i,j}$ belong to $(x_1,...,x_m)$. Then
\begin{small}
 \[
\xymatrix{
   R^{s}\ar[r]^{(a_{i,j}^{[k]})}& R^{t}\ar[r]&I^{[k]}\ar[r]&0
}
\]
\end{small}
is a minimal presentation of $I^{[k]}$ and all the entries
$a_{i,j}^{[k]}$ belong to $(x_1^k,...,x_m^k)\subset J$. Clearly,
this does not depend on the characteristic of the base field
$\Bbbk$.
\endproof
\
\\
\\
Next, the following lemma gives a sufficient condition for the
\emph{minimality} of the reduced horse resolution of $R/IJ$ (with
respect to $I$). Recall that $m_i(A)$ (resp. $M_i(A)$) is the
minimum (resp. maximum) shift in $i$-th homological degree in the
minimal graded free resolution of $A$.

\begin{lemma}\label{L2:HC}
 Let $(I,J)$ be a \emph{conjoined} pair of ideals in R. If
\[
 m_{k+1}(I)> M_1(I)+M_k (J)
\]
for $1\leq k \leq \text{pdim}(I) -1$, then the \emph{reduced
horseshoe resolution} of $R/IJ$ with respect to $I$ is minimal.
\end{lemma}
\emph{Proof.} Let
\begin{small}
\[
 \xymatrix{
0\ar[r]&G_m \ar[r]^{d_m^{''}}&\dots\ar[r]& G_2\ar[r]^{d_2^{'}} &
G_1\ar[r]^{d_1^{'}}&R\ar[r]^{\epsilon_{0}^{'}}& R/J\ar[r]&0
                      }
                      \]
\end{small}
be the minimal graded free resolution of $R/J$ and let
\begin{small}
\[
 \xymatrix{
0\ar[r]&F_n \ar[r]^{d_n^{''}}&\dots\ar[r]& F_2\ar[r]^{d_2^{''}} &
F_1\ar[r]^{d_1^{''}} & R\ar[r]^{\epsilon_{0}^{''}}& R/I \ar[r]&0
                      }
                      \]
\end{small}
be the minimal graded free resolution of $R/I$. Since the pair
$(I,J)$ is \emph{conjoined}, there is a minimal presentation of
$I$,
\begin{small}
 \[
\xymatrix{
   R^{s}\ar[r]^{\phi}& R^{t}\ar[r]^{\psi}&I\ar[r]&0,
}
\]
\end{small}
such that all the entries of the matrix that represents $\phi$
belong to $J$. Clearly, $F_1 \cong R^{t}$. Now, tensoring the
above exact sequence with $R/J$ yields
\[
  I/IJ \cong R^t \otimes R/J \cong F_1 \otimes R/J.
\]
Therefore, we get the minimal graded free resolution of $I/IJ$,
\[
 \xymatrix{
0\ar[r]&F_1 \otimes G_m \ar[r]&\dots\ar[r]& F_1 \otimes G_2\ar[r]
& F_1 \otimes G_1\ar[r]&F_1\ar[r]& I/IJ\ar[r]&0
                      }
                      \]
Next, note that for $j>1$ the maps $\lambda_j$ that appear in the
horseshoe lemma,
\[
  \lambda_j:F_j \to F_1 \otimes G_j
\]
are graded of degree zero. Since
\[
 m_j(I)> M_1(I)+M_{j-1} (J),
\]
we see that the degree of every basis element in $F_j$ is larger
than the degree of any basis element in $F_1 \otimes G_j$.
Therefore, the matrix that represents $\lambda_j$ does not have
any nonzero constant entry. Finally,
\[
  |G(IJ)|=|G(I)||G(J)|=dim(F_1)dim(G_1)=dim(F_1 \otimes G_1),
\]
and therefore the \emph{reduced horseshoe resolution} of $R/IJ$
with respect to $(I,J)$ is minimal.
\endproof

\begin{ex}
 Consider the \emph{conjoined} pair $(I,J)$ of the ideals  $I=(a^2 ,b^2, c^2 )$ and
 $J=(a,b,c)$ as in example \ref{EX:HC}. The Betti diagrams for the minimal resolutions of $R/(a,b,c)$ and
$R/(a^2,b^2,c^2)$ are
\[
\begin{array}{ccccc}
  total: & 1 & 3 & 3 & 1 \\
  0: & 1 & 3 & 3 & 1\\
\end{array}
\qquad
\begin{array}{ccccc}
total: & 1 & 3 & 3 & 1 \\
  0: & 1 & . & . & . \\
  1: & . & \textcolor{red}{3} & . & .\\
  2: & . & . & \textcolor{blue}{3} & .\\
  3: & . & . & . &\textcolor{green}{1}\\
  \end{array}
\]
In MACAULAY 2 we observe that the Betti diagram for the minimal
resolution of $R/J$ is
\[
\begin{array}{ccccc}
total: & 1 & 9 & 12 & 4 \\
  0: & 1 & . & . & . \\
  1: & . & . & .& . \\
  2: & . & \boxed{9=\textcolor{red}{3}\cdot 3} & \boxed{12=\textcolor{red}{3}\cdot 3 + \textcolor{blue}{3}}  & \boxed{3=\textcolor{red}{3} \cdot 1} \\
  3: & . & . & . & \textcolor{green}{1}
\end{array}
\]

The minimal graded free resolution of $R/I$ is of the form
\begin{small}
\[
\xymatrix{ 0\ar[r]&
R(-6)\ar[r]&R^3(-4)\ar[r]&R^3(-2)\ar[r]&R\ar[r]&R/I\ar[r]&0}.
\]
\end{small}
while the minimal graded free resolution of $R/J$ is of the form
\begin{small}
\[
\xymatrix{ 0\ar[r]&
R(-3)\ar[r]&R^3(-2)\ar[r]&R^3(-1)\ar[r]&R\ar[r]&R/J\ar[r]&0}.
\]
\end{small}
If we tensor this with $R^3(-2)$, we get the minimal graded free
resolution of $I/IJ$,
\begin{small}
\[
\xymatrix{ 0\ar[r]&
R^3(-5)\ar[r]&R^9(-4)\ar[r]&R^9(-3)\ar[r]&R^3(-2)\ar[r]&I/IJ\ar[r]&0}.
\]
\end{small}
Accordingly the \emph{reduced horseshoe resolution} of $R/IJ$ with
respect to $I$ is
\begin{small}
\[
\xymatrix{ 0\ar[r]& R(-6)\oplus R^3(-5)\ar[r]&R^{12}(-4)
\ar[r]&R^9(-3)\ar[r]&R\ar[r]&R/IJ\ar[r]&0},
\]
\end{small}
which is minimal.

\end{ex}

\section{Special ($p$-Borel) ideals}\label{SPB:S}

Let $p$ be a prime number and let $s,t$ be positive integers with
$p$-adic representations
\[
s=\sum a_i p^i \quad \text{and} \quad t=\sum b_i p^i
\]
with $0\leq a_i,b_i <p$. Then, we define the following order
$\prec_p$:
\begin{align}\notag
   s \prec_p t \quad &\iff {t \choose s} \ne 0\, \text{mod}\,p\\\notag
                     &\iff {a_i \leq b_i}\quad \text{for all} \,
                     i.\\\notag
\end{align}
If $x_j^t$ is the highest power of $x_j$ that divides a monomial
$\mathbf{m}$, we write $x_j^{t} || \mathbf{m}$.

\begin{defi} (see, e.g., \cite{Eisenbud} or \cite{Par}). $I$ is $p$-\emph{Borel}
if for every minimal generator $\mathbf{m}$ of $I$ and every
            $x_j$ such that $x_j^{t} || \mathbf{m}$, then
\[
        \left(\frac{x_{i}}{x_{j}}\right)^{s} \mathbf{m}
    \]
is in $I$ for all $i <j$ and $s \prec_p t$.\\
\end{defi}

Let $S=\{\mathbf{m}_1,\mathbf{m}_2,...,\mathbf{m}_r \}$ be a
finite set of monomials. If $I$ is the smallest $p$-Borel fixed
ideal such that $S$ is a subset of $G(I)$, then we say that $I$ is
generated by $\mathbf{m}_1,\mathbf{m}_2,...,\mathbf{m}_r$ in the
Borel sense and we write
\[
    I=<\mathbf{m}_1,\mathbf{m}_2,...,\mathbf{m}_r>.
\]
In particular, if $S=\{\mathbf{m}\}$, then $I$ is called
\textbf{principal $p$-Borel} and we write $I=<\mathbf{m}>$.

\begin{ex}. Let $R=\Bbbk[a,b,c]$ with $char (\Bbbk)=2$. Then the ideal
\[
I=(a^3,a^2b,ab^2,b^3,a^2c,b^2c,ac^2,bc^2)
\]
is a 2-Borel fixed ideal, minimally generated (in the Borel sense)
by $b^2c$ and $bc^2$; that is,
\[
I=<b^2c,bc^2>.
\]
 The ideal
\begin{align}\notag
    J&=(a,b,c)(a^2,b^2,c^2)\\\notag
    &=(a^3,a^2b,ab^2,b^3,a^2c,b^2c,ac^2,bc^2,c^3)\\\notag
      &=(I,c^3)\\\notag
      &=<c^3>\notag
\end{align}
is a principal 2-Borel fixed ideal.
\end{ex}

The first class of $p$-Borel ideals in $R=\Bbbk[x_1,...,x_n]$ that
were studied were the ones of the form
\[
  A=<x_n^{\mu}>,
\]
where $\mu$ is a positive integer, i.e. the Cohen-Macaulay
$p$-Borel fixed ideals (see, e.g. \cite{Ekki}, \cite{HerPop} and
\cite{Par}). The basic structure theory of principal $p$-Borel
ideals was developed in \cite{Par}, where it was proved that if
$\mu_k=\sum_{k,i} \mu_{ki}p^i$, where $0\leq \mu_{ki}\leq p-1$,
then
\[
  <x_1^{\mu_1}\cdots  x_n^{\mu_n}>=\prod_{k,i}
   (x_1^{p^i},...,x_k^{p^i})^{\mu_{ki}}.
\]
Products like the ones in the above structure of principal Borel
that depend on certain values of the $\mu_k$'s were studied in
\cite{Ene} and \cite{Wel}. These results and the observation of
the Betti diagrams of several $p$-Borel ideals in MACAULAY 2 led
us to study the following ideals.

\begin{defi} A \emph{special} ideal over $R$ is an ideal of the form
\[
    I=\prod_{j=1}^{s} I_{j}^{[p_j]},
\]
where
\[
I_j=(x_1,x_2,...,x_{\ell_j})^{a_j} \quad \text{and} \quad
I_j^{[p_j]}=(x_1^{p_j},x_2^{p_j},...,x_{\ell_j}^{p_j})^{a_j}
 \]
with
\[
    n= \ell_1 \geq \ell_2\geq ...\geq \ell_s \geq 1,
\]
and
\[
    0\leq a_j <\dfrac{p_{j+1}}{p_j}
    \]
with the numbers $\dfrac{p_{j+1}}{p_j}$ being integers $>1$ for
all $j=1,...,s$. In particular, if $p_j=p^{r_j}$ for $j=1,...,s$
for some prime number $p$ and some integers $r_s>...>r_2>r_1\geq
0$, we call it \emph{special} $p$-Borel ideal.
\end{defi}

 A \emph{special} $p$-Borel ideal is Borel
fixed if $char(k)=p$. Every $p$-Borel Cohen-Macaulay ideal is
\emph{special}, but as it is clear, not every principal $p$-Borel
ideal is \emph{special}.

\begin{ex}
Let $R=\Bbbk[a,b,c]$  ($char (\Bbbk)=2$). Then the 2-Borel ideal
$(a,b,c)^2 (a^4,b^4)$ is special, but not principal, since
\begin{align}\notag
  (a,b,c)^2 (a^4,b^4)&=(a^6,a^5b,a^4b^2,a^2b^4,ab^5,b^6,a^5c,a^4bc,ab^4c,b^5c,a^4c^2,b^4c^2)\\\notag
                     &= <b^5c,b^4c^2>.\notag
\end{align}
\end{ex}

The main result in this section is the following
\begin{thm}\label{P1:SPB}
The reduced horseshoe resolution of a \emph{special} ideal
$I=\prod_{j=1}^{s} I_{j}^{[p_j]}$ with respect to the ordered
$s$-tuple $\left(I_{s}^{[p_s]}, ...,I_{1}^{[p_1]}\right)$ is
minimal.
\end{thm}

In order to prove this theorem, we will introduce some notation
and prove lemma \ref{L1:SPB} and lemma \ref{L2:SPB} first. The
ideals $I_j$ are Borel fixed and their minimal graded free
 resolution is of the form
\begin{small}
\[
 \xymatrix{
           0\ar[r]& R^{\beta_{\ell_j, j}}(-d_{\ell_j,j}) \ar[r]&\dots\ar[r]& R^{\beta_{1,j}}(-d_{1,j})\ar[r]&R\ar[r]& R/I_{j} \ar[r]&0
                      }
                      \]
\end{small}
where $d_{i,j}=a_j+i-1$ for $i=1,2,...,\ell_j$ (see, e.g.
\cite{Elia}). Accordingly, the minimal graded free resolution of
$R/I_j^{[p_j]}$ is
\begin{small}
\[
 \xymatrix{
           0\ar[r]& R^{\beta_{\ell_j,j}}(-c_{\ell_j,j})\ar[r]^(0.7){\psi_{\ell_j,j}}&\dots\ar[r]& R^{\beta_{1,j}}(-c_{1,j})\ar[r]^(0.7){\psi_{1,j}}&R\ar[r]&R/I_j^{[p_j]}\ar[r]&0
                      }
                      \]
                      \end{small}
where $c_{i,j}=p_j d_{i,j}$ for $i=1,2,...,\ell_j$. This does not
depend on the characteristic of the base field $k$. Consider the
following free $R$-modules $F_{i,k}$,
\[
F_{i,1}=R^{\beta_{i,1}}(-(a_1+i-1) p_1)\qquad (1\leq i \leq
\ell_1)\]
 and for $1 < k \leq s$.
\[ F_{i,k}=
\begin{cases}
 \left(R^{\beta_{1,k}}(-a_k p_k) \otimes F_{i,k-1}\right) \oplus R^{\beta_{i,k-1}}(-(a_k+i-1) p_k) , &\text{if $1\leq i \leq \ell_k$;}\\
 R^{\beta_{1,k}}(-a_k p_k) \otimes F_{i,k-1},&\text{otherwise}.
 \end{cases}
\]
The $R$-module $F_{i,k}$ is the free module that appears in
homological degree $i$ in the minimal graded free resolution of
the ideal $J_k= I_{k-1}^{p_{k-1}} ...I_2^{p_2} I_1^{p_1}$ $(1\leq
k \leq s)$. The degrees of the basis elements of $F_{i,k}$, i.e.
the shifts in the minimal free resolution of $J_k$ are the
elements of the sets $S_{i,k}$, where for $1 \leq k \leq s$ we set
\[
S_{1,k}=\{a_1 p_1+a_2 p_2+...+a_k p_k \},
\]
and for $2\leq i \leq \ell_1=n $
\[
S_{i,k}=\{(a_1+i-1)p_1 + a_2 p_2+...+a_k p_k,
          (a_2+i-1)p_2 + a_3 p_3+...+a_k p_k,...,(a_k+i-1)p_k \}
\]
for $2\leq i \leq \ell_k -1$, and $S_{i,k}=S_{i,k-1}+\{a_k p_k\}$
otherwise.

Now we prove the following lemma.

 \begin{lemma}\label{L1:SPB}. We have
\begin{itemize}
\item[(a)] $(a_k +1)p_k > a_1 p_1 + a_2 p_2+...+a_k p_k$, for
$1\leq k \leq s$.
 \item[(b)] $(a_k +i-1)p_k > (a_{k-1}+i-2) p_{k-1}+a_k p_k$ for
  $2 \leq i\leq \ell_1$ and $k\geq 2$.
  \item[(c)] For fixed $i,k$ with $1 \leq k \leq s$ and $2\leq i \leq
  \ell_k -1$, the maximum element of $S_{i,k}$ is $c_{i,k}=(a_k+i-1)p_k$.
  \item[(d)] For fixed $k$ with $1 \leq k \leq s$, $\text{max}\{S_{i,k}-i|1\leq i \leq
\ell_1 -1\}$ is equal to $\text{max}\{S_{i,k}-i|\ell_k \leq i \leq
\ell_1 -1\}$, which is equal to the maximum of the elements
  \begin{align}\notag
(a_1+\ell_1-1)p_1+a_2 p_2+...+a_k p_k -\ell_1,\\\notag
 (a_2+\ell_2-1)p_2+a_3 p_3+...+a_k p_k-\ell_2,\\\notag
 ...,\\\notag
 (a_k+\ell_k-1)p_k-\ell_k.\notag
  \end{align}
 \end{itemize}
\end{lemma}
\textbf{\emph{Proof}}.
\begin{itemize}
    \item[(a)] By induction on $k$. For $k=1$, the inequality is
    trivially true. So assume that it is true for $k-1$ for some
    $k>1$. Then, by the induction hypothesis and since $a_{k-1} +1
    \leq \dfrac{p_k}{p_{k-1}}$, we get
    \begin{align}\notag
 a_1 p_1 + a_2 p_2+...+a_{k-1} p_{k-1}+a_k p_k  &< (a_{k-1}+1)p_{k-1}+a_k p_k\\\notag
                                                &\leq p_k+a_k p_k\\\notag
                                                &= (a_k +1) p_k,\notag
    \end{align}
    as desired.
    \item[(b)] Since $a_{k-1} p_{k-1}< p_k$ and  $p_k > p_{k-1}$, for $k>1$, it follows that
    for $1 \leq i\leq \ell_1$ we have
    \begin{align}\notag
     (a_k +i-1)p_k&=a_k p_k + (i-2)p_k +p_k\\\notag
                    &> a_k p_k  +(i-2)p_{k-1}+a_{k-1} p_{k-1}\\\notag
                    &= a_k p_k+(a_{k-1} +i-2)p_{k-1},\notag
    \end{align}
    as desired.
    \item[(c)] For $k=1$ this is clearly true, so assume that $k\geq 2$.
    Then it suffices to show that
    \[
(a_j+i-1)p_j + a_{j+1} p_{j+1}+...+a_k p_k \leq
(a_{j+1}+i-1)p_{j+1} + a_{j+2} p_{j+2}+...+a_k p_k
    \]
    for $1\leq j \leq k-1$. The above inequality is equivalent to
    \[
    (a_j+i-1)p_j \leq (i-1)p_{j+1},
    \]
    which is true, since
    \begin{align}\notag
     (a_j +i-1)p_j&=(a_j +1)p_j + (i-2)p_j\\\notag
                    &\leq p_{j+1} + (i-2)p_j\\\notag
                    &\leq p_{j+1} + (i-2)p_{j+1}\\\notag
                    &=(i-1)p_{j+1}.\notag
    \end{align}
    \item[(d)] First, the fact that $\text{max}\{S_{i,k}-i|1\leq i \leq
\ell_1 \}$ is equal to $\text{max}\{S_{i,k}-i|\ell_k \leq i \leq
\ell_1 \}$ follows from part (c) and
\[
  (a_k+i-1)p_k-i \geq (a_k+i-2)p_k-(i-1)
\]
for $i>2$ along with
\[
    (a_k+1)p_k-2 \geq (a_1 p_1 + a_2 p_2+...+a_{k-1} p_{k-1}+a_k p_k)-1,
\]
which is true from part (a). Now, for $k=1$, our claim is clearly
true, so assume that for
    some $k\geq 2$, the maximum element of $\text{max}\{S_{i,k-1}-i|1\leq i \leq
\ell_1 \}$ is equal to the maximum of the elements
  \begin{align}\notag
(a_1+\ell_1-1)p_1+a_2 p_2+...+a_{k-1} p_{k-1} -\ell_1 ,\\\notag
 (a_2+\ell_2-1)p_2+a_3 p_3+...+a_{k-1}p_{k-1}-\ell_2,\\\notag
 ...,\\\notag
 (a_{k-1}+\ell_{k-1}-1)p_{k-1}-\ell_{k-1}.\notag
  \end{align}
  From the definition of the sets $S_{i,k}$, we see that
  $\text{max}\{S_{i,k}-i|1\leq i \leq \ell_1
\}$ is equal to the maximum of
\[
 \text{max}\{S_{i,k}-i|1\leq i \leq \ell_k \},\]
 and
 \[
 \text{max}\{S_{i,k-1}-i|\ell_k \leq i \leq \ell_1\}+a_k p_k
\]
Now since
\[
 \text{max}\{S_{i,k}-i|1\leq i \leq \ell_k
 \}=(a_k+\ell_k-1)p_k-\ell_k,
 \]
and
\[
 \text{max}\{S_{i,k-1}-i|\ell_k \leq i \leq \ell_1\}+a_k p_k=\text{max}\{S_{i,k-1}-i|1 \leq i \leq \ell_1\}+a_k
 p_k,
\]
from the induction hypothesis, it follows that
$\text{max}\{S_{i,k}-i|1\leq i \leq \ell_1 \}$ is equal to the
maximum of the elements
\begin{align}\notag
(a_1+\ell_1-1)p_1+a_2 p_2+...+a_{k-1} p_{k-1}+a_k p_k
-\ell_1,\\\notag
 (a_2+\ell_2-1)p_2+a_3 p_3+...+a_{k-1} p_{k-1}+a_k
 p_k-\ell_2,\\\notag
 ...,\\\notag
 (a_{k-1}+\ell_{k-1}-1)p_{k-1}+ a_k p_k-\ell_{k-1}\\\notag
 (a_k+\ell_k-1)p_k-\ell_k.\notag
  \end{align}
  The proof is now complete.\endproof
\end{itemize}

In order to prove theorem \ref{P1:SPB} we need the following
lemma.

\begin{lemma} \label{L2:SPB}. Let $J_k= I_{k-1}^{p_{k-1}} ...I_2^{p_2} I_1^{p_1}$ $(1\leq k
\leq s)$. Then
 \begin{itemize}
 \item[(a)] All elements of $G(J_k)$ are of equal degree ($1\leq k \leq s$).
 \item[(b)] $\left(I_{k}^{[p_k]},J_{k-1} \right)$ is a \emph{conjoined} pair of
ideals for $1\leq k \leq s$. \item[(c)] $m_{j}(I_{k}^{[p_k]})>
M_1(I_{k}^{[p_k]})+M_{j-1} (J_{k-1})$ for $2\leq j \leq \ell_k$
and $2\leq k \leq s$.
\end{itemize}
\end{lemma}
\emph{Proof.}
\begin{itemize}
 \item[(a)] This follows by an easy induction and an argument as the
 one in the proof of Lemma \ref{L1:HC}.
 \item[(b)] This follows from the above part, Lemma \ref{L1:HC}
 and part (a) of Lemma \ref{L1:SPB}.
\item[(c)] Note that $m_{j}(I_{k}^{[p_k]})=(a_k+j-1)p_k$,
$M_1(I_{k}^{[p_k]})=a_k p_k$, $M_{j-1}
(J_{k-1})=(a_{k-1}+j-2)p_{k-1}$. Thus, the required inequality
becomes
\[
(a_k +j-1)p_k > a_k p_k+(a_{k-1}+j-2) p_{k-1},
\]
which is part (b) of Lemma \ref{L1:SPB}.\endproof
\end{itemize}
\
\\
\\

\textbf{\emph{Proof of Theorem \ref{P1:SPB}}}.  Let $2\leq k \leq
s$ and assume by induction that the minimal graded free resolution
of $R/J_{k-1}$ has been obtained already
\begin{small}
\[
 \xymatrix{
 0\ar[r]&F_{n,k-1} \ar[r]^(0.6){\phi_{n,k-1}}&...\ar[r]&F_{1,k-1}\ar[r]^(0.6){\phi_{1,k-1}}&R\ar[r] &R/J_{k-1} \ar[r]&0\\
                       }\]
\end{small}
where the degrees of the basis elements of the free $R-$modules
$F_{i,k-1}$ are the elements of $S_{i,k-1}$, for $1\leq i \leq n$.
We  know that the minimal graded free resolution of
$R/I_k^{[p_k]}$ is of the form
\begin{small}
\[
 \xymatrix{
           0\ar[r]& R^{\beta_{n,k}}(-c_{\ell_k,k})\ar[r]^(0.6){\psi_{\ell_k,k}}&\dots\ar[r]& R^{\beta_{1,k}}(-c_{1,k})\ar[r]^(0.6){\psi_{\ell_1,k}}&R\ar[r]&R/I_k^{[p_k]}\ar[r]&0
                      }
                      \]
\end{small}
From Lemma \ref{L2:SPB} part (b) above, we see that the pair
$(I_k^{[p_k]},J_{k-1})$ is \emph{conjoined}. Then Lemma
\ref{L2:SPB} part (c) together with Lemma \ref{L2:HC} implies that
the \emph{reduced horseshoe resolution} of $R/ I_k^{[p_k]}J_{k-1}$
is minimal and is of the form
\begin{small}
\[
 \xymatrix{
 0\ar[r]&F_{n,k} \ar[r]^(0.6){\phi_{n,k}}&...\ar[r]&F_{1,k}\ar[r]^(0.6){\phi_{1,k}}&R\ar[r] &R/J_{k} \ar[r]&0\\
                       }\]
\end{small}
\endproof

Since our techniques do not depend on the characteristic of the
base field, we obtain the following.
\begin{cor}
The graded betti numbers of a \emph{special} $p$-Borel ideal do
not depend on the characteristic of the base field $\Bbbk$.
\end{cor}

Also, we recover Pardue's regularity formula.
 \begin{cor}
Let $I=\prod_{j=1}^{s} I_{j}^{[p^{r_j}]}$ be a \emph{special}
$p$-Borel ideal. Then the regularity $reg(I)$ of $I$ is the
maximum of
\begin{align}\notag
    a_1 p^{r_1}+a_2 p^{r_2}+...+a_s p^{r_s}+(p^{r_1}-1)(\ell_1-1),\\\notag
    a_2 p^{r_2}+...+a_s p^{r_s}+(p^{r_2}-1)(\ell_2-1),\\\notag
 ...,\\\notag
 a_s p^{r_s}+(p^{r_s}-1)(\ell_s-1).\notag
  \end{align}
 \end{cor}
\textbf{\emph{Proof}}. This is immediate from the definition of
regularity combined with Lemma \ref{L1:SPB} part (d), the above
proposition and the fact that $reg(I)=reg(S/I)+1$.
\endproof
\section{Cellular resolutions} \label{S:CR}

In $char(\Bbbk)=0$, it is known that the Eliahou-Kervaire
resolution of a Borel fixed ideal $I$ is a CW-resolution (see
\cite{Batz}). If $I$ is generated in one degree, then another
minimal free resolution of $I$ can be supported on a polyhedral
cell complex (see \cite{Achilleas}). In $char(\Bbbk)=p$, it has
been proved in \cite{Ekki} that the \emph{minimal} free resolution
of a Cohen-Macaulay $p$-Borel fixed ideal is a CW-resolution and
in \cite{Wel} that the \emph{minimal} free resolution of a special
$p$-Borel fixed ideal is also a CW-resolution.

Here we construct a polyhedral cell complex that supports the
minimal free resolution of some special ideals. Our main result is
the following

\begin{prop} \label{T:CR} There exists a  polyhedral cell
complex that supports a minimal free resolution of a special ideal
$I$ of the form
\[
I=(x_1^{p_1},x_2^{p_1},...,x_{n}^{p_1})^{a_1}(x_1^{p_2},x_2^{p_2},...,x_{n}^{p_2})\cdots
(x_1^{p_s},x_2^{p_s},...,x_{n}^{p_s})
\]
\end{prop}

Before we prove this proposition, we consider a generalized
permutohedron ideal. Set $d:=a_1$ and recall that $dp_1 <p_2$. Let
$\mathbf{u}=(dp_1,p_2,0,...,0)$ be in $\mathbb{N}^n$. By permuting
the coordinates of $\mathbf{u}$, we obtain $n(n-1)$ points in
$\mathbb{N}^n$ constituting the vertices of an $(n-1)$-dimensional
\emph{generalized permutohedron} $\Pi (\mathbf{u})$ (see also,
\cite{Miller}). We label the vertices of $\Pi(\mathbf{u})$ by the
monomial generators of
\[
    K(\mathbf{u}):=(x_i^{dp_1}x_j^{p_2} \big| 1\leq i, j\leq n, i \ne
    j),
    \]
in a natural way and then we label an arbitrary face $F$ of
$\Pi(\mathbf{u})$ as usual, that is, by the lcm of the monomial
labels on all vertices in $F$. The inequality description of $\Pi
(\mathbf{u})$ is
\[
 \Pi(\mathbf{u})=\{(v_1,...,v_n) \in \mathbb{R}^n \big| \sum_{j=1}^n v_j=dp_1+p_2\, \text{and} \, 0\leq v_i \leq p_2 \, \text{for} \,
 i=1,2,...,n\},
\]
i.e. $\Pi (\mathbf{u})$ is the intersection of the $(n-1)$-simplex
\[
\Delta_{n-1}(x_1^{dp_1+p_2},\dots,x_n^{dp_1+p_2})=\{(v_1,...,v_n)
\in \mathbb{R}_{\geq 0}^n \big| \sum_{j=1}^n v_j=dp_1+p_2\}.
\]
 with the $n$ half
spaces $\{(v_1,...,v_n) \in \mathbb{R}^n \big| v_i \leq p_2\}$ (
$i=1,2,...n$). Since $K(\mathbf{u})_{\preceq \mathbf{b}}$ is empty
or contractible for all $\mathbf{b} \in \mathbb{Z}^n$,  $\Pi
(\mathbf{u})$ supports a free resolution of $K(\mathbf{u})$. It is
easy to see that that this resolution is minimal, since any two
comparable faces of the same degree coincide (see \cite{DS}).
Thus, we have proved the following.

\begin{lemma}\label{L1:CR}
The polyhedral cell complex $\Pi (\mathbf{u})$  supports a minimal
free resolution of $K(\mathbf{u})$.
\end{lemma}

We need the following lemma from \cite{Achilleas}.

\begin{lemma} \label{L2:CR} Let $I$ and $J$ be two monomial ideals in
$R$ such that $G(I+J)=G(I) \cup G(J)$ set-theoretically. Suppose
that
\begin{itemize}
 \item[(i)] $X$ and $Y$ are regular cell complexes in some
 $\mathbb{R}^{N}$ that support a (minimal) free resolution for
            $I$ and $J$, respectively, and
 \item[(ii)] $X\cap Y$ is a regular cell complex that supports a
            (minimal) free resolution for $I \cap J$.
            \end{itemize}
Then $X\cup Y$ supports a (minimal) free resolution for $I+J$.
\end{lemma}

 \textbf{Proof of Proposition \ref{T:CR}}. It suffices to consider the case $s=2$;
the general case is similar. Consider all monomial generators of
$(x_1^{p_1},...,x_n^{p_1})^d$ that are not divisible by $x_j$ and
denote the ideal they generate by $K_j$ ($1\leq j\leq n$). A
minimal free resolution of $K_j$ is supported on the
$(n-2)$-dimensional complex
$P_{j}:=P_{d}(x_1^{p_1},...,\widehat{x_j^{p_1}},...,x_n^{p_1})$.
Multiplying all vertices of $P_j$ by $x_j^{p_2}$, we obtain a
polyhedral cell complex $Q_j$ that supports a minimal free
resolution of $x_j^{p_2} K_j$.

Let $\sigma \subset [n]=\{1,2,...,n\}$. Replacing the face of $\Pi
(\mathbf{u})$ that lies on the hyperplane $\{(v_1,...,v_n) \in
\mathbb{R}^n \big| v_j =p_2\}$ by $Q_j$ for $j \in \sigma$ gives
us a polyhedral cell complex $\Pi_{\sigma}$ that supports a
minimal free resolution of the ideal
\[
  K(\mathbf{u}) +\sum_{j \in \sigma} x_j^{p_2} K_j .
\]

The intersection of $\Pi_{\sigma}$ with $x_j^{p_2}
P_{a_1}(x_1,...,x_n)$ is $Q_j$ for $j \in \sigma$. Applying lemma
\ref{L2:CR}, we glue all these complexes to obtain a polyhedral
cell complex that supports a minimal free resolution of
\[
 K(\mathbf{u}) +\sum_{j \in \sigma} x_j^{p_2}(x_1^{p_1},x_2^{p_1},...,x_n^{p_1})^{d}.
\]
In particular, when $\sigma=[n]$, we obtain a polyhedral cell
complex that supports a minimal free resolution of
\[
  (x_1^{p_1},x_2^{p_1},...,x_n^{p_1})^{d}(x_1^{p_2},x_2^{p_2},...,x_n^{p_2}).
\]
\endproof

\begin{rem} It follows from the above proposition that there exists a polyhedral
cell complex that supports a minimal free resolution of any
Cohen-Macaulay 2-Borel fixed ideal.
\end{rem}

\begin{ex} The polyhedral
cell that supports a minimal free resolution of
$(a,b,c)(a^2,b^2,c^2)$ is
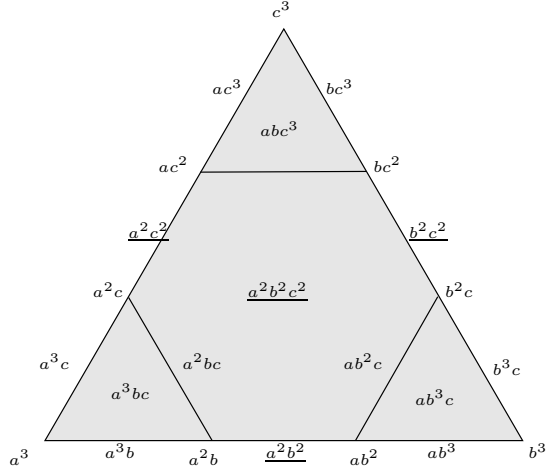
\begin{figure}[htbp]
\begin{center}
 \input{celfig3.pstex_t}
\end{center}
\end{figure}
\end{ex}

\section{Mapping cone}\label{S:mc}

By applying the iterated mapping cone technique (see \cite{Xara},
\cite{HerTan}, \cite{PS}) as in the case of Borel fixed ideals in
characteristic zero, we do not always obtain a minimal resolution
in characteristic $p$. The smallest example that we found in
characteristic two using MACAULAY 2 \cite{Mike} is the following
one in three variables.

\begin{ex}\label{e:mc}
Let $R=\Bbbk[a,b,c]$ ($char (\Bbbk)=2$),
\[
I=(a^3,a^2b,ab^2,b^3,a^2c,b^2c,ac^2,bc^2),
\]
and
\[
I'=(a^3,a^2b,ab^2,b^3,a^2c,b^2c,ac^2).
\]
Starting with the ideal $(a^3)$ and adding the monomial generators
of $I'$ one at a time in the order that appears above, the
iterated mapping cone gives us a \emph{minimal} free resolution of
$I'$. However, if $f$ is the map from the resolution of
$R/(I':bc^2)$ to the resolution of $R/I'$ induced by
multiplication by $bc^2$, then the mapping cone of $f$ does not
give us a minimal free resolution of $R/I$. This is clear from the
following Betti diagrams of $I'$ and $I$, since
$\beta_{2,5}(I')=1$, while $\beta_{2,5}(I)=0$.
\[
 {\begin{array}{ccccc}
   \text{total}: & 1 & 7 & 9 & 3  \\
  0 :& 1 & . & . & . \\
  1 :& . & . & . & .\\
  2 :& . & 7 & 8 & 2 \\
  3 :& . & . &1 & 1 \\
\end{array}} \quad
{\begin{array}{ccccc}
  \text{total}: & 1 & 8 & 10 & 3  \\
  0 :& 1 & . & . & . \\
  1 :& . & . & . & .\\
  2 :& . & 8 & 10 & 2 \\
  3 :& . & . & . & 1 \\
\end{array}}\]
\end{ex}

As the proof of the following proposition shows, the ordering of
the monomial generators of the above ideal is not important.

\begin{prop}
    There exists a $p$-Borel fixed ideal $I$
such that for any ordering $m_1\succ m_2 \succ...\succ m_r$ of its
minimal generators, there is some $i$ with $2\leq i \leq r$, such
that the mapping cone of multiplication by $m_i$ from a minimal
resolution of $R/((m_1,...,m_{i-1}):m_i)$ to a minimal resolution
of $R/(m_1,...,m_{i-1})$ is not a minimal free resolution of
$R/(m_1,...,m_i)$.
\end{prop}
\textbf{Proof}. Let $R=\Bbbk[a,b,c]$ with $char (\Bbbk)=2$ and
consider the ideal
\[
I=(a^3,a^2b,ab^2,b^3,a^2c,b^2c,ac^2,bc^2),
\]
The minimal cellular resolution of $I$ consists of two triangles
with vertices in $\{a^3,a^2b,a^2c\}$ and $\{ab^2,b^3,b^2c\}$, and
the hexagon $\Pi(1,2,0)$ with vertices in $\{a^2b,ab^2,a^2c,
b^2c,ac^2,bc^2\}$.

\begin{figure}[htbp]
\begin{center}
 \input{celfig4.pstex_t}
\end{center}
\end{figure}

It suffices to check whether we can get the above hexagon by an
iterated method, i.e. by adding one monomial at a time in a
suitable order and by considering each time the minimal cellular
resolution of the corresponding ideal that we get.

This is impossible. Indeed, let $I(S)$ be the ideal generated by
be a 5-element subset $S$ of the vertices of the 6-gon
$\{a^2b,ab^2,a^2c, b^2c,ac^2,bc^2\}$. By considering cases, it is
easy to see that the total betti number $\beta_3(I(S))$ is
non-zero. That is, the minimal free resolution of $I(S)$ is
supported on a 2-dimensional polygon, i.e. there is an edge that
connects two vertices, which are not connected in  $\Pi(1,2,0)$.
For example, if $S= \{a^2b,ab^2,a^2c, b^2c,ac^2\}$, we see that
there is an edge between $ac^2$ and $b^2c$, which is denoted with
a dashed edge in the above figure. This means that when we add
$bc^2$ to $S$, we must erase that edge in order to obtain
$\Pi(1,2,0)$, which supports a minimal free resolution of the
ideal generated by the new set $S\cup \{bc^2\}$.
\endproof

\begin{rem}
This means that there is a lifting map $\lambda_i$ in the mapping
cone of the map $f$ from example \ref{e:mc}  with a non-zero
constant entry in the matrix that represents it. Thus, there is a
multigraded free module $R(-\mathbf{e})$ that appears in the
$i$-th homological degree in both resolutions of $R/I'$ and
$R/(I':bc^2)$. If we cancel (by a change of basis) the two copies
of $R(-\mathbf{e})$ that appear in the mapping cone, we obtain a
minimal free resolution of $R/I$.

However, this does not mean that every time we have a copy of
$R(-\mathbf{e})$ in the same homological degree, we could obtain a
minimal free resolution by cancelling it. One of the smallest
examples in characteristic two that we found using MACAULAY 2
\cite{Mike} is the following one in five variables.

\begin{ex}
Let $R=\Bbbk[a,b,c,d,e]$  with $char (\Bbbk)=2$ and let $B$ be the
2-Borel fixed ideal
\[
  B=<ace^2,b^2e^2,bcd^2,bc^2d>.
\]
Then
\[
  B:(bce^2)=(b,a,d^2,cd,c^2)
\]
and the 2-Borel fixed ideal $A=(B,bce^2)=<bc^2d,bce^2>$ has 30
generators. The Betti diagrams of $B$ and $A$ are\\
\[
 {\begin{array}{ccccccc}
   \text{total}: & 1 & 29 & 78 & 83 & 41 & 8\\
  0 :& 1 & . & . & . &.&.\\
  1 :& . & . & . & .&.&.\\
  2 :& . & .& . & . &.&.\\
  3 :& . & 29 & 56 & 34 & 6 & .\\
  4 :& . & . & 22 & 48 & 30 & 5\\
  5 :& . & . & . & 1 & 5 & 3\\
\end{array}} \quad
{\begin{array}{ccccccc}
   \text{total}: & 1 & 30 & 83 & 91 & 46 & 9\\
  0 :& 1 & . & . & . &.&.\\
  1 :& . & . & . & .&.&.\\
  2 :& . & .& . & . &.&.\\
  3 :& . & 30 & 58 & 35 & 6 & .\\
  4 :& . & . & 25 & 56 & 36 & 6\\
  5 :& . & . & . & . & 4 & 3\\
\end{array}}
\]
We note that one copy of each of $R(-(1,2,1,0,2))$,
$R(-(0,2,2,2,2))$ and  $R(-(1,2,2,2,2))$ appears in homological
degrees 2, 4 and 5, respectively, in the resolutions of $R/B$ and
$R/(B:bce^2)$, but the two copies of $R(-(1,2,1,0,2))$ that appear
in the mapping cone cannot be cancelled.
\end{ex}
\end{rem}

\textbf{Acknowledgements}: I would like to express my thanks to
Mike Stillman for his support and for the numerous discussions we
have had on this project. Also, I would like to thank Irena Peeva
for suggesting the study of $p$-Borel fixed ideals as an
interesting problem.

Address:

 Cornell University, Department of Mathematics,
 Ithaca, NY 14853,

 asin@math.cornell.edu

\end{document}